\documentclass[12pt]{amsart}
\usepackage{amscd,amssymb,amsmath}
\newtheorem{theorem}{Theorem}[section]
\newtheorem{lemma}[theorem]{Lemma}
\newtheorem{proposition}[theorem]{Proposition}
\newtheorem{corollary}[theorem]{Corollary}

\theoremstyle{definition}

\newtheorem{example}[theorem]{Example}

\theoremstyle{remark}
\newtheorem{remark}{Remark}[section]

\numberwithin{equation}{section}
\begin{document}
\title{Twisted conjugacy classes in nilpotent groups}
\author{Daciberg Gon\c calves}
\address{Dept. de Matem\'atica - IME - USP, Caixa Postal 66.281 - CEP 05311-
970,
S\~ao Paulo - SP, Brasil; FAX: 55-11-30916183}
\email{dlgoncal@ime.usp.br}
\author{Peter Wong}
\address{Department of Mathematics, Bates College, Lewiston,
ME 04240, U.S.A.; FAX: 1-207-7868331}
\email{pwong@bates.edu}
\thanks{This work was initiated during the first author's visit to Bates College April 8 - 30, 2006.
The paper was completed during the second author's visits to IME-USP July 9 - 20, October 20 - 23, 2006.
The visit of the first author was supported by FAPESP, Projeto Tem\'atico
Topologia Alg\'ebrica, Geom\'etrica e Differencial-2000/05385-8 and the visits of the second author were supported by the National Science Foundation OISE-0334814.}
\begin{abstract}
A group is said to have the $R_\infty$ property if every automorphism has an infinite number of twisted conjugacy classes. We study the question whether $G$ has the $R_\infty$ property when $G$ is a finitely generated torsion-free nilpotent group. As a consequence, we show that for every positive integer $n\ge 5$, there is a compact nilmanifold of dimension $n$ on which every homeomorphism is isotopic to a fixed point free homeomorphism. As a by-product, we give a purely group theoretic proof that the free group  on two generators has the $R_\infty$ property. The $R_{\infty}$ property for virtually abelian and for $\mathcal C$-nilpotent groups are also discussed.
\end{abstract}
\date{\today}
\dedicatory{Dedicated to Professor A. Dold and to Professor E. Fadell}
\maketitle

\newcommand{\af}{\alpha}
\newcommand{\et}{\eta}
\newcommand{\ga}{\gamma}
\newcommand{\ta}{\tau}
\newcommand{\ph}{\varphi}

\newcommand{\bt}{\beta}

\newcommand{\lb}{\lambda}

\newcommand{\wh}{\widehat}

\newcommand{\sg}{\sigma}

\newcommand{\om}{\omega}

\newcommand{\cH}{\mathcal H}

\newcommand{\cF}{\mathcal F}

\newcommand{\N}{\mathcal N}

\newcommand{\R}{\mathcal R}

\newcommand{\Ga}{\Gamma}

\newcommand{\cc}{\mathcal C}

\newcommand{\bea} {\begin{eqnarray*}}

\newcommand{\beq} {\begin{equation}}

\newcommand{\bey} {\begin{eqnarray}}

\newcommand{\eea} {\end{eqnarray*}}

\newcommand{\eeq} {\end{equation}}

\newcommand{\eey} {\end{eqnarray}}

\newcommand{\ovl}{\overline}

\newcommand{\vv}{\vspace{4mm}}

\newcommand{\lra}{\longrightarrow}


\bibliographystyle{plain}

\section{Introduction}

A natural question with a long history in group theory is whether an infinite group must have an infinite number of conjugacy classes of elements. In 1949, Higman-Neumann-Neumann constructed an infinitely generated group with only a finite number of conjugacy classes. Subsequently, S. Ivanov constructed one such example where the group is finitely generated. More recently, D. Osin gave an example of a finitely generated infinite group in which all non-trivial elements are conjugate. More generally, given a group endomorphism $\varphi: \pi \to \pi$, one considers the $\varphi$-twisted conjugacy classes. Equivalently, $\pi$ acts on $\pi$ via $\sigma \cdot \alpha \mapsto \sigma \alpha \varphi(\sigma)^{-1}$ and the Reidemeister number $R(\varphi)$, the cardinality of the set of orbits of this action, is the same as the number of $\varphi$-twisted conjugacy classes. Clearly, $R({\rm id}_{\pi})$ is simply the number of conjugacy classes. From the point of view of Nielsen-Wecken fixed point theory, $R(\varphi)$ is precisely the number of fixed point classes of any map $f:X\to X$ with induced homomorphism $f_{\#}=\varphi$ on $\pi_1(X)=\pi$. For many spaces $X$ (e.g. $X$ is a compact manifold of dimension at least three), the Nielsen number $N(f)$, which is the principal object of study, coincides with the minimal number $MF[f]$ of fixed points of maps within the homotopy class of $f$. The Nielsen number is always bounded above by the Reidemeister number $R(f)=R(\varphi)$ and for a large class of spaces, $N(f)=0$ or $N(f)=R(f)$. While the computation of $N(f)$ is in general very difficult, the problem of determining the finiteness of $R(\varphi)$ is more tractable.
In 1985, D. Anosov \cite{anosov} and independently E. Fadell and S. Husseini \cite{ed-suf} showed that for any selfmap $f:N\to N$ of a compact nilmanifold, $|L(f)|=N(f)$ where $L(f)$ denotes the Lefschetz number of $f$. Using the approach of \cite{ed-suf}, this was later strengthened in \cite{felshtyn-hill-wong} where it was shown in particular that $N(f)>0$ iff $R(f)<\infty$. In fact, for selfmaps $f$ of a nilmanifold, either $N(f)=0$ or $N(f)=R(f)$.

In 1994, A. Fel'shtyn and R. Hill \cite{fel-hill} conjectured that for a finitely generated group $\pi$ of exponential growth, if $\varphi:\pi \to \pi$ is injective then $R(\varphi)=\infty$. Using techniques from geometric group theory, G. Levitt and M. Lustig \cite{levitt-lustig} showed that if $\pi$ is finitely generated torsion-free non-elementary Gromov hyperbolic then every automorphism $\varphi \in {\rm Aut}(\pi)$ must have $R(\varphi)=\infty$. This was subsequently proven in \cite{fel:1} without the torsion-free assumption.
We say that a group {\it $G$ has the
$R_{\infty}$ property for automorphisms}, in short  {\it $G$ has the
$R_{\infty}$ property}, if for every automorphism $\varphi:G \to G$ we have $R(\varphi)=\infty$. It was shown in \cite{daci-peter} that the Fel'shtyn-Hill conjecture does not hold in general. Moreover, non-elementary polycyclic groups of exponential growth that are not Gromov hyperbolic without the $R_{\infty}$ property were constructed. Since then, many examples of groups with $R_{\infty}$ have been discovered (see e.g., \cite{fel-daci}, \cite{fgw}, \cite{felshtyn-leonov-troitsky}, \cite{daci-peter4}, \cite{levitt}, \cite{TW1}, and \cite{TW2}). (For connections between the study of Reidemeister classes and other areas, see e.g. \cite{felshtyn-troitsky}.) In particular, groups that are quasi-isometric to generalized Baumslag-Solitar groups  \cite{TW2} have the $R_{\infty}$ property. In these examples, the groups are all of exponential growth. Since a result of M. Gromov states that a finitely generated group is of polynomial growth iff it is virtually nilpotent, it is natural to ask if one can find (virtually) nilpotent groups with the $R_{\infty}$ property. The main objective of this paper is to construct finitely generated groups of polynomial growth that have the $R_{\infty}$ property. This shows that the $R_{\infty}$ property does {\it not} depend on the growth type of the group as suggested by the original conjecture of Fel'shtyn and Hill.

It is easy to see that finitely generated abelian groups do not have the $R_{\infty}$ property. Therefore, we will first explore the $R_{\infty}$ property for virtually abelian and for nilpotent groups.

The basic algebraic techniques used in the present paper for showing $R(\varphi)=\infty$ is the relationship among the Reidemeister numbers of groups homomorphisms of a short exact sequence. In general, given a commutative diagram of groups and homomorphisms
\begin{equation*}
\begin{CD}
    A    @>{\eta}>>  B  \\
    @V{\psi}VV  @VV{\ph}V   \\
    A    @>{\eta}>>  B
\end{CD}
\end{equation*}
the homomorphism $\eta$ induces a function $\hat {\eta}:\mathcal R(\psi) \to \mathcal R(\ph)$ where $\mathcal R(\alpha)$ denotes the set of $\alpha$-twisted conjugacy classes.
Some of the basic facts that will be used throughout this paper are given in the following lemma. For more general results, see \cite{daci-peter} and \cite{wong}.

\begin{lemma}\label{R-facts}
Given an endomorphism $\psi:G\to G$ of a finitely generated torsion-free abelian group $G$, $Fix \psi=\{1\}$ iff $R(\psi)<\infty$.
Consider the following commutative diagram
\begin{equation*}\label{general-Reid}
\begin{CD}
    1 @>>> A    @>>>  B @>>>    C @>>> 1 \\
    @.     @V{\ph'}VV  @V{\ph}VV   @V{\ovl \ph}VV @.\\
    1 @>>> A    @>>>  B @>>>    C @>>> 1
 \end{CD}
\end{equation*}
where the rows are short exact sequences of groups and the vertical arrows are group endomorphisms.

(1) If $R(\ovl \ph)=\infty$ then $R(\ph)=\infty$.

(2) If $R(\ovl \ph)<\infty, |Fix \ovl \ph|<\infty$, and $R(\ph')=\infty$ then $R(\ph)=\infty$.

(3) If the short exact sequence is a central extension then $R(\varphi)=R(\varphi')R(\bar {\varphi})$.
\end{lemma}
\begin{proof}
Suppose $G$ is finitely generated torsion-free abelian. Then $G=\mathbb Z^k$ for some positive integer $k$. For any $\varphi: G\to G$, $R(\varphi)=\# Coker (1-\varphi)=|\det (1-\varphi)|$ so that $R(\varphi)<\infty$ iff $\varphi$ does not have $1$ as an eigenvalue iff $\varphi(x)=x$ has only trivial solution, i.e., $Fix \varphi=1$.

The homomorphism $p:B\to C$ induces a function $\hat p:\mathcal R(\varphi) \to \mathcal R(\ovl{\varphi})$ given by $\hat p([\alpha]_B)=[p(\alpha)]_C$. Since $p$ is surjective, so is $\hat p$. Thus, $(1)$ follows. Similarly, $i:A\to B$ induces a function $\hat i:\mathcal R(\varphi') \to \mathcal R(\varphi)$. Since the sequence is exact, it follows that $\hat i(\mathcal R(\varphi'))=\hat p^{-1}([1]_C)$. The subgroup $Fix \ovl{\varphi}$ acts on $\mathcal R(\varphi')$ via $\bar {\theta}\cdot [\alpha]_A=[\theta \alpha \varphi(\theta)^{-1}]_A$ where $\theta \in B$ and $\bar {\theta}\in Fix \ovl{\varphi}$. Thus, two classes $[\alpha]_A$ and $[\beta]_A$ are in the same $Fix \ovl{\varphi}$-orbit iff $i(\alpha)$ and $i(\beta)$ are in the same Reidemeister class, i.e., $[i(\alpha)]_B=[i(\beta)]_B$. Now, $(2)$ follows immediately. Finally, if the extension is central, $\hat p^{-1}([\bar \alpha]_C)$ is independent of $\bar \alpha$ so that $R(\varphi)=k\cdot R(\varphi')$ and $k$ is the number of distinct Reidemeister classes of $\ovl{\varphi}$ in $C$. In other words, $k=R(\ovl{\varphi})$ and thus $(3)$ follows.
\end{proof}

In \cite{dyer}, J. Dyer constructed some nilpotent Lie algebras which have the property that
every automorphism is unipotent. This implies that the induced homomorphism on the abelianization is the identity. This should imply that every automorphism of the corresponding nilpotent group has the property that the induced automorphism on the abelianization is the identity. An immediate consequence of this is the fact that every automorphism $\varphi$ of such nilpotent group has an infinite number of $\varphi$-twisted conjugacy classes and thus the group has the
$R_{\infty}$ property. The following example is the group theoretic analog of the Lie algebra theoretic example constructed by Dyer \cite{dyer}.

\begin{example}\label{Dyer}
Let $G$ be the quotient of the  free nilpotent group $F_2/\Gamma_7(F_2)$ on two generators $x,y$ by the normal  closure of the words
$Y_1^{-1}Y_3$, $[Y_2,x]$, $[Y_1,y]U_4^{-1}$, $[Y_1,x][Y_2,y]^{-3}U_4^{-1}$. Here $B=[x,y]$, $Y_1=[[B,y],y]$, $Y_2=[[B,y],x]$, $Y_3=[[B,x],x]$, $z_1=[Y_1,y]$, $z_2=[Y_1,x]$, $z_3=[Y_2,y]$, $U_4=[z_3,x]$.  Then given any automorphism
$\varphi$ of $G$ the induced homomorphism  on the abelianization of $G$ is the identity. Therefore, this group has the $R_{\infty}$ property. It has nilpotency class 6 and Hirsch length $\geq 7$.
\end{example}

\begin{remark} One can construct as in \cite{bryant-papistas} a finitely generated torsion-free nilpotent group $G$ so that ${\rm Aut}(G)$, modulo the kernel of the action of ${\rm Aut}(G)$ on the abelianization $G^{Ab}$, is a prescribed closed subgroup of $GL(2,\mathbb Z)$.
\end{remark}

This example has prompted us to determine which finitely generated torsion-free nilpotent groups can have the $R_{\infty}$ property.
In general, the $R_{\infty}$ property has important implications in the fixed point
theory of Jiang-type spaces since $R(f)=\infty$ implies that $N(f)=0$ which in turn implies in most cases that $f$ is deformable to be fixed point free. Recall that a space is said to be of Jiang-type if $L(f)=0 \Rightarrow N(f)=0$ and $L(f)\ne 0 \Rightarrow N(f)=R(f)$. Since nilmanifolds are known to be of Jiang-type, for each  torsion-free nilpotent group  which has the $R_{\infty}$ property there corresponds a nilmanifold with the property that every homeomorphism is homotopic to a fixed point free map. Such examples cannot appear in dimension 3 but for every $n>3$ there exists an $n$-dimensional nilmanifold such that every homeomorphism is homotopic to be fixed point free; and if $n\ge 5$ then the homotopy can be made to be isotopy. As a by-product of our investigation, we are able to give an algebraic proof that automorphisms of free groups on two generators have infinite Reidemeister number.

This paper is organized into seven sections. In section 2, extension of the $R_{\infty}$ property to automorphisms of virtually abelian groups is discussed. In particular, we show that the fundamental group of the Klein bottle (Theorem \ref{Klein}), which is a finitely generated torsion-free virtually abelian group, has the $R_{\infty}$ property. We also construct for every $n \ge 2$ a finitely generated torsion-free virtually abelian group of Hirsch length $n$ with the desired $R_{\infty}$ property (Theorem \ref{Klein-Zn}). In section 3 we show that $F_2/\Gamma_k$ has the $R_{\infty}$ property for $k\ge 9$ (Theorem \ref{free-k-R}) where $F_2$ is the free group on two generators and $\Gamma_k=\Gamma_k(F_2)$ is the $k$-th term of the lower central series of $F_2$. In sections 4 and 5, examples of finitely generated torsion-free nilpotent groups with the $R_{\infty}$ property are presented, according to the nilpotency class and the Hirsch length. In particular, we show that for each $n>3$, there is a nilpotent group of Hirsch length $n$ with the $R_{\infty}$ property. Furthermore, nilpotent groups of nilpotency class 2 that have the $R_{\infty}$ property are also constructed.
In the section 6, we turn to topological applications. We show that there exist nilmanifolds of dimension $n$ for each $n\ge 5$ on which every homeomorphism is isotopic to a fixed point free homeomorphism (Theorem \ref{fpf-maps-on-nilmanifolds}). In the last section, we consider a generalization of nilpotent groups, namely $\mathcal C$-nilpotent groups where $\mathcal C$ denotes the class of finite groups. These groups provide further examples of groups with the $R_{\infty}$ property. For any $n\ge 7$, we make use of the results from section 5 to construct a compact $\mathcal C$-nilpotent manifold of dimension $n$ for which every homeomorphism is isotopic to a fixed point free homeomorphism (Theorem \ref{fpf-maps-on-C-nilpotent-spaces}).

\section{virtually abelian groups}

In this section, we study the $R_{\infty}$ property for virtually abelian groups. These groups are finite extensions of abelian groups and they appear as fundamental groups of almost flat Riemannian manifolds.

The simplest torsion-free virtually abelian group with the $R_{\infty}$ property is the fundamental group of the Klein bottle $K$ which is the total space of a $S^1$-bundle over $S^1$. On the other hand, $K$ is finitely covered by the $2$-torus $T^2$ so that
$\pi_1(K)$ has $\mathbb Z^2$ as a finite index subgroup. Since $\mathbb Z^2$ does not have the $R_\infty$ property, it follows that in general the $R_\infty$ property is not even invariant under commensurability. Note that $\pi_1(K)$ is non-elementary non-Gromov hyperbolic and has polynomial growth.
Furthermore, we will show that the groups $\pi_1(K)\times \mathbb Z^n$ for
every $n\geq 0$, which are infra-abelian having $\mathbb Z^{n+2}$ as a subgroup of
index 2, have the $R_{\infty}$ property.

We will start by analyzing the group $\mathbb Z\rtimes \mathbb Z$ where the action is non-trivial. This group is the fundamental group of the Klein bottle.
The group of automorphisms of $\mathbb Z\rtimes \mathbb Z$ is known as folklore. Since we cannot find a precise reference for this fact, we state it in the following lemma and we sketch a proof for completeness.

\begin {lemma} The group ${\rm Aut}(\mathbb Z\rtimes \mathbb Z)$ where $\mathbb Z\rtimes \mathbb Z$ is as above is isomorphic to
$(\mathbb Z\bigoplus \mathbb Z_2)\rtimes \mathbb Z_2$ where the action  $\theta:\mathbb Z_2 \to {\rm Aut}(\mathbb Z\bigoplus \mathbb Z_2)$ sends the generator to the automorphism  $(r,\epsilon) \to (-r,\epsilon)$. Furthermore, the inner automorphisms are isomorphic to  the subgroup $(2\mathbb Z\bigoplus 0)\rtimes \mathbb Z_2$ and the quotient ${\rm Out}(G)$ is isomorphic to $\mathbb Z_2\oplus \mathbb Z_2$.
\end{lemma}
\begin{proof}
We use the presentation $\langle \alpha,\beta|\alpha \beta \alpha \beta^{-1} \rangle$ for the fundamental group of the Klein bottle. Given an endomorphism $\kappa$, write $\kappa(\alpha)=\alpha^{\epsilon}\beta^s$ and $\kappa(\beta)=\alpha^{r}\beta^{\delta}$. Now,
\begin{equation}\label{trivial}
\begin{aligned}
\kappa(\alpha \beta \alpha \beta^{-1})&=\kappa(\alpha)\kappa (\beta) \kappa(\alpha) \kappa(\beta^{-1}) \\
                                      &=\alpha^{\epsilon}\beta^s \alpha^{r}\beta^{\delta} \alpha^{\epsilon}\beta^s (\alpha^{r}\beta^{\delta})^{-1} \\
                                      &=\alpha^q\beta^{s+\delta+s-\delta}=\alpha^q\beta^{2s}
\end{aligned}
\end{equation}
for some integer $q$ that depends on the values of $\epsilon, s, r$, and $\delta$. Since $\alpha \beta \alpha \beta^{-1}=1$, the expression \eqref{trivial} yields the trivial element. Thus, $2s=0$ and hence $s=0$ so that $\kappa(\alpha)=\alpha^{\epsilon}$. It follows that
the automorphisms of this group are of the form $\alpha \mapsto \alpha^{\epsilon}, \beta \mapsto \alpha^r \beta^{\delta}$. Moreover, $\epsilon, \delta \in \{1,-1\}$. To see this, we note that given any two automorphisms $\alpha \mapsto \alpha^{\epsilon_1}, \beta \mapsto \alpha^{r_1} \beta^{\delta_1}$ and
$\alpha \mapsto \alpha^{\epsilon_2}, \beta \mapsto \alpha^{r_2} \beta^{\delta_2}$, the composite is an automorphism  such that $\alpha \mapsto \alpha^{\epsilon_1\epsilon_2}$
and the $\beta$-exponent of the image of $\beta$ is $\delta_1\delta_2$. Consider the automorphisms
\begin{equation*}
\begin{aligned}
&\varphi_1: \qquad \alpha \mapsto \alpha ; \quad \beta \mapsto \alpha \beta \\
&\varphi_2: \qquad \alpha \mapsto \alpha ; \quad \beta \mapsto \beta^{-1} \\
&\varphi_3: \qquad \alpha \mapsto \alpha^{-1} ; \quad \beta \mapsto \beta
\end{aligned}
\end{equation*}
Writing $\mathbb Z_2=\{\pm 1\}$, the generators for $(\mathbb Z \oplus \mathbb Z_2)\rtimes \mathbb Z_2$ are $\eta_1=(1,1,1), \eta_2=(0,-1,1)$ and $\eta_3=(0,1,-1)$. The homomorphism $\varphi_i \mapsto \eta_i$ for $i=1,2,3$ is the desired isomorphism.
\end{proof}

Now we will show that $\mathbb Z \rtimes \mathbb Z$ has the $R_{\infty}$ property.

\begin{theorem}\label{Klein}
For any automorphism $\varphi$ of  $\mathbb Z\rtimes \mathbb Z$, we have  $R(\varphi)=\infty$.
\end{theorem}
\begin{proof}  Denote by $x,y$ the generators of the group $\mathbb Z\rtimes \mathbb Z$, where
$x$ is a generator of the first copy of $\mathbb Z$ and $y$ is a generator of the
second copy of $\mathbb Z$. In this group, which is the   fundamental group of the
Klein bottle, we have the relation  $xyxy^{-1}=1$. The automorphisms of
$\mathbb Z\rtimes \mathbb Z$ from their  description given in the proof of Lemma 2.1   can be divided into four cases as follows.\\
 a) $x\mapsto x; y\mapsto x^ry~$ \ \ \ \ \ \        b) $x\mapsto x; y \mapsto x^ry^{-1}$ \\
 c) $x \mapsto x^{-1}; y \mapsto  x^ry$ \ \ \ \   d) $x\mapsto x^{-1}; y \mapsto
x^ry^{-1}$

\noindent
where $r\in \mathbb Z$.

Cases a) and c) are treated as follows. These automorphisms leave invariant the subgroup generated by $x$.
So every such automorphism induces a homomorphism of the short exact sequence
$$0\to \mathbb Z\to \mathbb Z\rtimes \mathbb Z \to \mathbb Z\to 0$$
and induces in the quotient the identity homomorphism of $\mathbb Z$ which has an infinite number
of conjugacy classes. So the result follows.

For Case b), we have the automorphism  $x\to x; y \to x^ry^{-1}$
which maps a generic element $x^my^k$ to $x^my^{-k}$ if $k$ is even and to
$x^{m+r}y^{-k}$ if $k$ is odd (these elements are obtained  by a straightforward
calculation using the relation $xyxy^{-1}=1$ on the group).
 So the elements of the group in the Reidemeister class of the element $x^i$
are of the form either
$x^my^{2n}x^iy^{2n}x^{-m}$ or $x^my^{2n+1}x^iy^{2n+1}x^{-m-r}$.
After simplifying the expression using the relation of the group,
these elements can be written as  $x^iy^{4n}$ or  $x^{-r-i}y^{4n+2}$, respectively. So
there are an infinite number of distinct Reidemeister classes since  $x^i$ is in the Reidemeister class
of $x^j$ if and only if $i=j$. This proves Case b).

For Case d), we have the automorphism  $x\to x^{-1}; y \to x^ry^{-1}$
which maps a generic element $x^my^k$ to $x^{-m}y^{-k}$ if $k$ is even and to
$x^{-m+r}y^{-k}$ if $k$ is odd.
 So the elements of the group in the Reidemeister class of the element $x^iy$
are the elements of the form either
$x^my^{2n}x^iyy^{2n}x^{m}$ or $x^my^{2n+1}x^iyy^{2n+1}x^{m-r}$.
Again, using the relation of the group,
these elements can be written as  $x^iy^{1+4n}$ or  $x^{r-i}y^{4n+3}$, respectively. So
there are an infinite number of distinct Reidemeister classes since  $x^iy$ and $x^jy$ are in the same Reidemeister class
if and only if $i=j$. This proves Case d).
\end{proof}

Next, we construct from the Klein bottle examples of Bieberbach groups with the $R_{\infty}$ property. First, the following result is crucial in our construction.

\begin{proposition}\label{simplest}
The group $\mathbb Z\rtimes \mathbb Z_2$ has the $R_{\infty}$ property.
\end{proposition}
\begin{proof} Let $\varphi:\mathbb Z\rtimes \mathbb Z_2\to \mathbb Z\rtimes \mathbb Z_2$ be an automorphism.
Since the only elements of infinite order are the elements of the form $t^r$
for $t\in \mathbb Z$ a generator and $r\ne 0$, the subgroup $\mathbb Z$ is characteristic.
Therefore, $\varphi|_{\mathbb Z}:\mathbb Z \to \mathbb Z$. Since ${\rm Aut}(\mathbb Z)=\{\pm 1\}$, $\varphi|_{\mathbb Z}$ is either the identity or multiplication by $-1$. In the former case, the assertion that $R(\varphi)=\infty$ follows from $(2)$ of Lemma \ref{R-facts}.

For the latter case, the Reidemeister class of $(t^l,1)$ contains at most two elements where $1$ is the non-trivial
element of $\mathbb Z_2$. To see this, suppose $(t^r,1)$ and $(t^s,1)$ are $\varphi$-conjugate, that is, $(t^s,1)=\alpha (t^r,1) \varphi(\alpha)^{-1}$ for some $\alpha=(t^j,\epsilon) \in \mathbb Z \rtimes \mathbb Z_2$ where $\epsilon=0,1$. Suppose $\varphi(t^0,1)=(t^n,1)$. Now, $\varphi(\alpha)=(t^{-j+\epsilon n},\epsilon)$ so that $\varphi(\alpha)^{-1}$ is equal to $(t^{-j+n},1)$ if $\epsilon =1$ or $(t^j,0)$ if $\epsilon = 0$. It follows that
$$\alpha (t^r,1) \varphi(\alpha)^{-1}=(t^j,1)(t^r,1)(t^{-j+n},1)=(t^{j-r},0)(t^{-j+n},1)=(t^{n-r},1) \qquad \text{when $\epsilon=1$}$$
and
$$\alpha (t^r,1) \varphi(\alpha)^{-1}=(t^j,0)(t^r,1)(t^{j},0)=(t^{j+r},1)(t^{j},0)=(t^{r},1) \qquad \text{when $\epsilon=0$}.$$ Since the set $\{(t^l,1)\}$ is infinite and the Reidemeister classes have finite length, it follows that $R(\varphi)=\infty$.
\end{proof}

The next result shows that for any positive integer $n\ge 2$, one can construct a finitely generated Bieberbach group of Hirsch length $n$ that has the $R_{\infty}$ property.

\begin{theorem}\label{Klein-Zn}
The group $\pi_1(K)\times \mathbb Z^n$ has the $R_{\infty}$ property for
all integer $n\geq 0$.
\end{theorem}
\begin{proof} Consider the presentation $\pi_1(K)=\langle a,b | abab^{-1}\rangle$.
The center of this group is the subgroup generated by $b^2$ and similarly the
center of $\pi_1(K)\times \mathbb Z^n$ is the subgroup $\langle b^2\rangle \times \mathbb Z^n$.
Now consider the short exact sequence
$$1\to \langle b^2\rangle \times \mathbb Z^n \to \pi_1(K)\times \mathbb Z^n \to \mathbb Z\rtimes \mathbb Z_2 \to 1.$$

Since the center is characteristic with  respect to automorphisms, any
automorphism of  $\pi_1(K)\times \mathbb Z^n$ is an automorphism of the short exact
sequence. Since the quotient  $\mathbb Z\rtimes \mathbb Z_2$ has the $R_{\infty}$
property by Proposition \ref{simplest}, the result follows from Lemma \ref{R-facts}.
\end{proof}

\begin{remark} This result in the case $n=0$ gives an alternate proof of Theorem \ref{Klein}.
\end{remark}

\begin{remark} Given a virtually nilpotent group $G$, we have a short exact sequence $1 \to N \to G \to F \to 1$ where $F$ is finite and $N$ is nilpotent. If $N$ is characteristic and has the  $R_{\infty}$ property then it follows easily from Lemma \ref{R-facts} that $G$ also has the $R_{\infty}$ property.
\end{remark}

\section{Commutators and free nilpotent groups}

Given $r$ a positive integer and  $k$ either an integer or $\infty$, consider the free nilpotent group $G(r,k)=F_r/\Gamma_{k+1}(F_r)$ of rank $r$ and nilpotency class $k$, where for $k=\infty$  we have $G(r,\infty)=F_r$ because
$F_r$ is residually nilpotent. For which values of $r$ and $k$(including $k=\infty$) does the group $G$ have the $R_\infty$ property? For the case $k=\infty$, the free group $F_r$ for $r\ge 2$ is a torsion-free non-elementary Gromov hyperbolic group thus it follows from \cite{levitt-lustig} that $F_r$ has the $R_{\infty}$ property. We should point out that the proof in \cite{levitt-lustig} is geometric. In this section,
we will show that $F_2/\Gamma_9(F_2)$ has the property. It seems that
for $i<9$ the group $F_2/\Gamma_i(F_2)$ does not have this property.
As a result, we will give a purely algebraic proof of
the fact that an automorphism of a free group of rank $2$ has an infinite number
of twisted conjugacy classes.

Let us denote $\Gamma_n(F_2)$ simply by $\Gamma_n$.
From Witt's Formulae (\cite{kms} Theorem 5.11), the ranks of the abelian groups $\Gamma_i/\Gamma_{i+1}$ are calculated. In particular, we have:\\
 $rk(\Gamma_2/\Gamma_{3})=1 $, $rk(\Gamma_3/\Gamma_{4})=2 $, $rk(\Gamma_4/
\Gamma_{5})=3$, $rk(\Gamma_5/\Gamma_{6})=6 $,
 $rk(\Gamma_6/\Gamma_{7})=9$, $rk(\Gamma_7/\Gamma_{8})=18 $, $rk
(\Gamma_8/\Gamma_{9})=30$.

\begin{lemma} Let $G(2,k)$ be the  free nilpotent group of rank $2$ on two generators $a$ and $b$,  nilpotency class $k$  and
 $\varphi:G(2,k) \to G(2,k)$  a homomorphism. Suppose $k\geq 3$. The induced homomorphism
$\tilde \varphi:\Gamma_2(G(2,k))/\Gamma_3(G(2,k))\to \Gamma_2(G
(2,k))/\Gamma_3(G(2,k))$ is multiplication by $\det (M)$ where $M$ is the matrix
of ${\varphi}^{Ab}:G(2,k)^{Ab} \to G(2,k)^{Ab}$ and $G(2,k)^{Ab}$ denotes the abelianization of $G(2,k)$.
\end{lemma}

\begin{proof} Since $k\ge 3$, it follows that $\Gamma_i(G(2,k))=\Gamma_i(F_2)/\Gamma_{k+1}(F_2)$ for $i\le k$. This means that $\Gamma_2(G(2,k))/\Gamma_3(G(2,k))\cong \Gamma_2/\Gamma_3$ and $G(2,k)^{Ab} \cong {F_2}^{Ab}$. Suppose the induced map ${\varphi}^{Ab}$ is given by the matrix $M=\left(\begin{smallmatrix}
                                                                \alpha & \beta \\
                                                                \gamma & \delta
                                                              \end{smallmatrix}\right)$.
Note that $\Gamma_2/\Gamma_3$ is generated by the single element $B\Gamma_3$ where $B=[a,b]=aba^{-1}b^{-1}$. Using classical calculus of commutators (see \cite{kms} Chapter 5), we have $Ba\Gamma_3=aB\Gamma_3$ and $Bb\Gamma_3=bB\Gamma_3$. Since $aba^{-1}=Bb$, it is straighforward to verify that $\tilde \varphi(B\Gamma_3)=B^k\Gamma_3$ where $k=\alpha \delta -\beta \gamma=\det(M)$.
\end{proof}

\begin{corollary}\label{det=1}
If $\varphi:G(2,k) \to G(2,k)$ is an automorphism such that $\det({\varphi}^{Ab})=1$
then $R(\varphi)=\infty$.
\end{corollary}
\begin{proof} We have a homomorphism of short exact sequence $0\to \Gamma_2/
\Gamma_3 \to G(2,k)/\Gamma_3 \to G(2,k)/\Gamma_2 \to 0$. Since the induced map on the
subgroup $\Gamma_2/\Gamma_3=\mathbb Z$ is multiplication by $1$, it follows from formula (2.2) of  \cite{daci-peter} that  $G(2,k)/\Gamma_3$ has the $R_{\infty}$ property. Since $\Gamma_3$ is characteristic in $G(2,k)$, the result follows for  $G(2,k)$.
\end{proof}

In order to study the $R_{\infty}$ property, it suffices, as a result of Corollary \ref{det=1}, to consider  automorphisms $\varphi$ whose induced automorphisms ${\varphi}^{Ab}$ have determinant $-1$.
While we have succeeded in determining the $R_\infty$ property for $G(2,k)$ for $k\geq 9$, we do not know how to extend our techniques to $G(r,k)$ for  $r>2$ and $k>2$. On the other hand, the case of the free abelian group $G(r,1)$ of rank $r$ is well-understood.

Next, we will construct a nontrivial element in $\Gamma_8/\Gamma_9\subset G(2,8)$ which is fixed by the restriction of any automorphism $\varphi:G(2,8)\to G(2,8)$.

\begin{proposition} If $\varphi:G(2,8)\to G(2,8)$ is an automorphism then $R(\varphi)=\infty$.
\end{proposition}
\begin{proof} By Corollary \ref{det=1}, the result holds if $\det(\varphi^{Ab})=1$. So let us
assume that  $\det(\varphi^{Ab})=-1$. Then consider the quotient $\Gamma_3/\Gamma_4$.
This group has rank 2 and has as generators $[B,x],[B,y]$ where $B=[x,y]$. A
straightforward calculation shows that the matrix of the automorphism induced
by $\varphi$ has determinant -1. Therefore the image of the element $[[B,x],[B,y]]
\in \Gamma_6/\Gamma_7$ is $([[B,x],[B,y]])^{-1}$. Call this element $w$. From
\cite{kms} Lemma 5.4 (pg 314), it follows that $w\ne 1$. Now consider the element
$w_1=[B,w]\in \Gamma_8/\Gamma_9$. From above, we have $\varphi'(w_1)=w_1$ where $\varphi'$ is the induced homomorphism on $\Gamma_8/\Gamma_9$ and
again from \cite{kms} Lemma 5.4 we conclude that $w_1\ne 1$. Then the automorphism $\varphi$ induces an automorphism $\hat \varphi$  of the following short exact sequences and hence a commutative diagram.
\begin{equation*}
\begin{CD}
0   @>>> \Gamma_8/\Gamma_9 @>>> F_2/\Gamma_9 @>>>     F_2/\Gamma_8               @>>> 1  \\
@.  @V{\varphi'}VV                 @V{\hat \varphi}VV          @V{\bar \varphi}VV               @.  \\
0   @>>> \Gamma_8/\Gamma_9 @>>> F_2/\Gamma_9 @>>>     F_2/\Gamma_8               @>>> 1
\end{CD}
\end{equation*}
If the
induced map on the quotient has infinite Reidemeister number then the result
follows. Otherwise, it will be finite. Note that $\varphi'(w_1)=w_1$ and $\Gamma_8/\Gamma_9$ is finitely generated torsion-free abelian. Applying Lemma \ref{R-facts}, we have $R(\varphi')$ is infinite and together with the fact that $R(\bar \varphi)$ is finite, we conclude that $R(\varphi)$ is infinite and the
result follows.
\end{proof}

Since the commutator subgroups $\Gamma_k$ are characteristic, we can use induction and the same argument as in the proof above to the short exact sequence $0\to \Gamma_k/\Gamma_{k+1} \to F_2/\Gamma_{k+1} \to F_2/\Gamma_k \to 1$ to show that $F_2/\Gamma_{k+1}$ has the desired property. Thus, we state the following result.

\begin{theorem}\label{free-k-R}
Let $G(2,k)$ be the free nilpotent group on $2$ generators of nilpotency class $k$. If $k\ge 9$ (including the case $k=\infty$) then for any automorphism $\varphi:G(2,k) \to G(2,k)$, we have $R(\varphi)=\infty$, i.e., $G(2,k)$ has the $R_\infty$ property.
\end{theorem}

\begin{remark} When $k=\infty$, $G(2,\infty)=F_2$. Thus, our proof of Theorem \ref{free-k-R} provides a purely algebraic proof of the fact that $F_2$ has the $R_{\infty}$ property as a consequence of a result proven in \cite{levitt-lustig} using techniques from geometric group theory.
\end{remark}

\section{$R_{\infty}$ and nilpotency class}

If $G$ is a finitely generated torsion-free abelian group, then $G$ does not have the $R_{\infty}$ property. Since abelian groups have nilpotency class $1$ and Example \ref{Dyer} gives an example where the nilpotency class is $6$, it is interesting to know if one can construct a nilpotent group of nilpotency class $\le 5$ that has the desired $R_{\infty}$ property. First, is there a finitely generated torsion-free nilpotent group of nilpotency class  $2$  which has
the $R_\infty$ property?

\begin{example}\label{nil-class2}
We now construct a group of nilpotency class 2 which has the
$R_{\infty}$ property. Consider the free nilpotent group $G(4,2)$ on $3$ generators $x,y,z,w$ of nilpotency class 2 and $\Gamma_2(G(4,2))/\Gamma_3(G(4,2))$ is a free abelian group of rank 6 in which the cosets defined by the commutators $[x,y],[x,z],[x,w]$,$[y,z],[y,w],[z,w]$ form a basis. Take the quotient by the subgroup generated by $[x,z],[x,w],[y,z]$. We will prove that this group has the $R_{\infty}$ property.

Let $\varphi:F_4/\Gamma_3\to F_4/\Gamma_3$ be an automorphism  such that the matrix on the abelianization is given by
$$
   M=\left[ \begin{array}{cccc}
        a_1 & a_2 & a_3 & a_4 \\
        b_1 & b_2 & b_3 & b_4 \\
        c_1 & c_2 & c_3 & c_4 \\
        d_1 & d_2 & d_3 & d_4
     \end{array} \right].
$$

The proof is divided into several steps. \\
Step 1. We will show that $b_1=d_1=b_3=d_3=0$. Suppose that  $b_3\ne 0$. The proof of the other cases is similar and we leave to the reader. Because the determinant of each of the $2\times 2$ matrices,

$$
   \left[ \begin{array}{cc}
        a_i & a_j \\
        b_i & b_j
     \end{array} \right],
\left[ \begin{array}{cc}
        b_i & b_j \\
        d_i & d_j
     \end{array} \right],
   \left[ \begin{array}{cc}
        c_i & c_j \\
        d_i & d_j
     \end{array} \right]
  $$

\noindent  for $(i,j)=(1,3),(2,3),(1,4)$ is zero
we have $(a_1,a_3)=\lambda(b_1,b_3)$,  $(a_2,a_3)=\lambda_1(b_2,b_3)$ so that
$\lambda=\lambda_1$. Similarly  $(d_1,d_2,d_3)=\alpha(b_1,b_2,b_3)$.
If $\alpha=0$ then $d_4\ne 0$ since the determinant of the matrix is nonzero.
But $c_1d_4=b_1d_4=0$ which implies that  $c_1=b_1=0$. So the first column of $M$ is
zero, a contradiction. So let $\alpha \ne 0$. Then similarly because
$\alpha b_3 \ne 0$, it follows that $(c_1,c_2, c_3)=\theta(b_1,b_2,b_3)$. So the rows of $M$ up to
the third column are proportional and the determinant is zero. This is a
contradiction and thus $b_3=0$.\\

Step 2. From now on, we will consider that  $b_1=d_1=b_3=d_3=0$. In this step we will compute the form of the matrix under the assumption that $a_3\ne 0$. Because  $b_1=d_1=b_3=d_3=0$ we must have $b_2a_3=d_2c_3=a_1b_4=c_1d_4=0$. Since
$a_3\ne 0$ implies  $b_2=0$, $b_4\ne 0$, and $a_1=0$, the matrix $M$ is of the form
$$
   M=\left[ \begin{array}{cccc}
        0 & a_2 & a_3 & a_4 \\
        0 & 0 & 0 & b_4 \\
        c_1 & c_2 & c_3 & c_4 \\
        0 & d_2 & 0 & d_4
     \end{array} \right].
$$

\noindent Now $\det (M)=b_4a_3c_1d_2$. Because this determinant is $\pm 1$, it follows that all the factors are either 1 or -1. Also following from the equations above, we have $c_3=d_4=0$. Thus
$\det (M-Id)=1-b_4d_2-c_1a_3+b_4d_2c_1a_3$. If  $\det (M-Id)=0$ then we have $R(\varphi^{Ab})=\infty$. This always happens except when  $b_4d_2=c_1a_3=-1$ in which case $\det (M-Id)=4$. In this situation, we are going to compute the matrix
$N$ of the induced homomorphism on $\Gamma_2/\Gamma_3$. This matrix is

  $$
N=\left[ \begin{array}{ccc}
        0 & a_2b_4 & a_3b_4  \\
        0 & -d_2b_4 & 0 \\
        c_1d_2 & -d_2c_4 & 0
     \end{array} \right].
$$
\noindent Using the fact that $d_2b_4=-1$, we have $\det (N-Id)=0$ since the second row becomes trivial and the result follows.

Step 3. Suppose now   that $a_3=0$. Then we have
 $$
   M=\left[ \begin{array}{cccc}
        a_1 & a_2 & 0 & a_4 \\
        0 & b_2 & 0 & b_4 \\
        c_1 & c_2 & c_3 & c_4 \\
        0 & d_2 & 0 & d_4
     \end{array} \right].
$$

\noindent So we have   $d_2c_3=a_1b_4=c_1d_4=0$. Since $c_3\ne 0$, it  folows that $d_2=0$. Since $\det (M)=a_1b_2c_3d_4$, it follows that all factors are $\pm 1$ and consequently $b_4=c_1=0$. Then $\det (M-Id)= (a_1-1)(b_2-1)(c_3-1)(d_4-1)$ and as  before if   $\det (M-Id)=0$ then the Reidemeister number is infinite. This always happens except when  $a_1=b_2=c_3=d_4=-1$ in which case $\det (M-Id)$ is 16. In this case we are going to compute the matrix
$N$ of the induced homomorphism on $\Gamma_2/\Gamma_3$. This matrix is

 $$
N=\left[ \begin{array}{ccc}
        a_1b_2 & -b_2a_4 & 0  \\
        0 & b_2c_4 & 0 \\
        0 & c_2d_4 & c_3d_4
     \end{array} \right].
$$
\noindent Using the fact that $a_1=b_2=-1$ we have $\det (N-Id)=0$ since the first column  becomes trivial and the result follows. This concludes the proof.
\end{example}

In contrast to the last example, we will next show that the {\it free} nilpotent groups of nilpotency class $2$ do not have the $R_{\infty}$ property.

\begin{example}\label{nil-class2}
We now show that the free nilpotent groups $G(r,2)$ does not have the $R_{\infty}$ property. It suffices for each $r$ to exhibit an automorphism
$\varphi:G(r,2) \to G(r,2)$ such that $R(\varphi)<\infty$.\\
 For $r=2$, consider the automorphism given by $\varphi(x)=x^2y$ and  $\varphi(y)=x^5y^2$. The matrix of the induced automorphism on the abelianization $\mathbb Z\oplus \mathbb Z$ has two eigenvalues $\lambda_1=2+5^{1/2}$, $\lambda_2=2-5^{1/2}$. It follows that $Fix(\varphi^{Ab})=\{1\}$ so that $R(\varphi^{Ab})<\infty$. The induced automorphism $\varphi'$ on $\Gamma_2/\Gamma_3=\mathbb Z$ is multiplication by $-1$ so that $R(\varphi')<\infty$. Since we have a central extension we have $R(\varphi)=R(\varphi^{Ab})R(\varphi')$ by Lemma \ref{R-facts} and the result follows.\\
Let $r\ge 3$. If $r$ is even,  define $\varphi:(\mathbb Z\oplus \mathbb Z)^{r/2} \to (\mathbb Z\oplus \mathbb Z)^{r/2}$ as a direct sum of the automorphism $\varphi:\mathbb Z\oplus \mathbb Z \to \mathbb Z\oplus \mathbb Z$ defined above. If $r$ is odd, define
$$\varphi:(\mathbb Z\oplus \mathbb Z)^{(r-1)/2}\oplus\mathbb Z \to (\mathbb Z\oplus \mathbb Z)^{(r-1)/2}\oplus\mathbb Z$$
as a direct sum of the automorphism $\varphi:\mathbb Z\oplus \mathbb Z \to \mathbb Z\oplus \mathbb Z$ defined above and the automorphism given by multiplication by -1 in the last coordinate. It is easy to see that the product of any two eigenvalues of $\varphi$ for either $r$ even or
$r$ odd  cannot be 1. Therefore
$\varphi':\Gamma_2/\Gamma_3\to \Gamma_2/\Gamma_3$ does not have $1$ as an eigenvalue so that $Fix \varphi'=\{1\}$ and the result follows from Lemma \ref{R-facts}.
\end{example}

\section{$R_{\infty}$ and Hirsch length}

There is only one finitely generated torsion-free nilpotent group of Hirsch length 2, namely
$\mathbb Z\oplus\mathbb Z$ and it does not have the $R_\infty$ property.
For which integers $n$, is there a finitely generated torsion-free nilpotent group of Hirsch length $n$ with the $R_\infty$ property?

\begin{example}\label{length-ex1}
Let us consider nilpotent
torsion-free groups of Hirsch length 3. They  are classified as follows. For each integer $r$, consider the group $N_r$ given by the following presentation $\langle a,b,c| [a,c]=1, [b,c]=1$ and  $[a,b]=c^r\rangle$. It is not hard to show that there is an automorphism
 $\varphi:N_r \to N_r$ defined by the map which sends  $a\to a^2b$ and $b \to a^5b^2$. The automorphism $\varphi$ induces a commutative diagram of automorphisms of the short exact sequence
$$0\to C\to N_r \to \mathbb Z \oplus \mathbb Z \to 0$$
\noindent where $C$ is the center of the group. Since the extension is central and the induced automorphisms $\varphi'$ and $\ovl {\varphi}$ have finite Reidemeister numbers, it follows from (4) of Lemma \ref{R-facts} that $R(\varphi)$ is finite.
\end{example}

\begin{example}\label{ex3}
We now construct a finitely generated nilpotent group $G$ of Hirsch length 4 such that
$\Gamma_4(G)=1$, $\Gamma_3(G)=\mathbb Z$, and $G$ has the $R_{\infty}$ property. This group is a quotient of
the nilpotent group $F_2/\Gamma_4(F_2)$ where $F_2$ is the free group on the
generators. The group $G$ has the property that
$G/\Gamma_3(G)=F_2/\Gamma_3(F_2)$.  The difference comes between
$\Gamma_3(G)/\Gamma_4(G)=\mathbb Z$ and   $\Gamma_3(F_2)/\Gamma_4(F_2)=\mathbb Z \oplus \mathbb Z$ (see
e.g. \cite{kms}). Now we describe $G$ in terms of generators and relations. Let
$x,y$ be the generators, $B=[x,y]$ and $w_1=[[x,y],x]$ such that
$[[x,y],y]=1$. Note that the last relation does not hold in the group $F_2/\Gamma_4(F_2)$.
Now given an automorphism of $G$, the matrix of the abelianization has a special
form since $[B,y]=1$. This matrix has the first column $a,b$ and second
column $0,d$. As a consequence, the automorphism induced on the abelianization
either has eigenvalue 1  or $\det ( \varphi)=1$. In the former case the result follows because the induced map on the abelianization has Reidemeister number infinite. In the latter case the result follows from an argument similar to Lemma 2.1. Let $H=F_2/\Gamma_4(F_2)$.  Then the group $G$ given by the presentation
$\langle  x, y|\Gamma_4,[B,y]\rangle$, which can be regarded as a quotient of $H$,  has the $R_\infty$ property.
\end{example}

\begin{remark} Example \ref{ex3} does not imply that $F_2$
has the $R_{\infty}$ property but it gives an example of a nilpotent group with the $R_{\infty}$ property that is simpler than $F_2/\Gamma_9(F_2)$ . Further, it has the least possible Hirsch length among all finitely generated torsion-free nilpotent groups having the $R_{\infty}$ property.
\end{remark}

\begin{example}\label{n>3}
Here we construct for every integer $\ell\geq 4$ a finitely generated torsion-free nilpotent group of Hirsch length $\ell$ with the $R_{\infty}$ property.  We will show that the  product of the group given in Example \ref{ex3} with $\mathbb Z^n$ for $n\geq 0$, has the desired property.

Now for any $n\ge 0$, let $G_n=G\times \mathbb Z^{n}$ where $G=\langle  x, y|\Gamma_4,[B,y]\rangle$ is the group constructed in Example \ref{ex3}. Then, $\Gamma_k(G_n)=\Gamma_k(G)$ so $G_n$ is of rank $n+2$, has  nilpotency class $3$, and its Hirsch length  is $n+4$. Moreover, $\Gamma_3(G_n)\subseteq Z(G_n)=Z(G) \times \mathbb Z^{n}$. Now, given $\varphi\in {\rm Aut}(G_n)$, we have the following commutative diagram.

\begin{equation}
\begin{CD}
    0 @>>> \Gamma_3(G_n) @>>> G_n           @>>> G/\Gamma_3(G)=G_n/\Gamma_3(G_n) @>>> 1   \\
    @.     @V{\varphi'}VV        @V{\varphi}VV          @V{\bar \varphi=\bar \zeta}VV          @.  \\
    0 @>>> \Gamma_3(G_n) @>>> G_n           @>>> G/\Gamma_3(G)=G_n/\Gamma_3(G_n) @>>> 1
\end{CD}
\end{equation}

Since $G$ has the $R_\infty$ property and the following commutative diagram

\begin{equation}
\begin{CD}
    0 @>>> \Gamma_3(G) @>>> G           @>>>      G/\Gamma_3(G)        @>>> 1   \\
    @.     @V{\zeta'}VV      @V{\zeta}VV          @V{\bar \zeta}VV          @.  \\
    0 @>>> \Gamma_3(G) @>>> G           @>>>      G/\Gamma_3(G)        @>>> 1
\end{CD}
\end{equation}
yields the product formula $R(\zeta)=R(\zeta')R(\bar \zeta)$ by $(4)$ of Lemma \ref{R-facts}, it follows that
$$R(\zeta)=\infty \Rightarrow R(\zeta')=\infty \quad \text{or} \quad R(\bar \zeta)=\infty.$$ Thus $R(\varphi)=\infty$ if $R(\bar \zeta)=\infty$. If $R(\zeta')=\infty$ then $R(\varphi')=\infty$ since $\Gamma_3(G_n)$ projects onto $\Gamma_3(G)$. Hence, $R(\varphi)=\infty$.\\
\end{example}

\section{Fixed point free homeomorphisms on nilmanifolds}

Recall that a compact nilmanifold is the coset space of a finite dimensional nilpotent Lie group by a closed cocompact subgroup. A classical result of A. Mal'cev asserts that such nilmanifolds are coset spaces of simply connected nilpotent Lie groups by uniform discrete subgroups. Furthermore, there is a one-to-one correspondence between finitely generated torsion-free nilpotent groups $\Gamma$ and compact nilmanifolds $M=K(\Gamma,1)$.

From Example \ref{n>3}, we can associate to each such finitely generated torsion-free nilpotent group a compact nilmanifold so that the Reidemeister number of each homeomorphism is infinite. As an application, we obtain the following result for fixed point free homeomorphisms on these nilmanifolds.

\begin{theorem}\label{fpf-maps-on-nilmanifolds}
For any $n\ge 4$, there exists a compact nilmanifold $M$ with $\dim M=n$ such that every homeomorphism $f:M\to M$ is homotopic to a fixed point free map. Moreover, if $n\ge 5$ then $f$ is isotopic to a fixed point free homeomorphism.
\end{theorem}
\begin{proof} For any $n\ge 4$, let $M$ be the compact $K(G_{n-4},1)$ nilmanifold where $G_k$ denotes the group as constructed in Example \ref{n>3}. Since the Hirsch length of $G_{n-4}$ is $n$, it follows that $\dim M=n$. If $f:M\to M$ is a homeomorphism, then $f_{\#}$ is an automorphism of $G_{n-4}$. Thus, $R(f_{\#})=\infty$. It follows from \cite{felshtyn-hill-wong} or \cite{go:nil2} that $N(f)=0$, that is, the Nielsen number vanishes. It follows from a classical theorem of Wecken that $f$ is deformable to be fixed point free. If $n\ge 5$, we invoke the isotopy Wecken theorem of Kelly \cite{kelly} to obtain the desired fixed point free homeomorphism in the isotopy class of $f$.
\end{proof}

\begin{remark} The isotopy Wecken theorem of Kelly also holds in dimension 2 \cite{jiang-guo} based upon Nielsen-Thurston classification of surface homeomorphisms and in dimension 3 \cite{jww} based upon techniques of $3$-manidols. It is however unknown whether it holds in dimension $4$ in general. In fact, it is not even known whether a homeomorphism of the $4$-torus $T^4$ with zero Lefschetz number (hence deformable to be fixed point free) can be isotopic to a fixed point free homeomorphism.

It is known \cite{MO} that the well-known Arnold conjecture holds for nilmanifolds which states that the minimal number of fixed points of Hamiltonian symplectomorphisms of a compact nilmanifold of dimension $2n$ is at least $2n+1$ (for more details on the Arnold conjecture, see e.g. \cite{R}). Based upon our results, any symplectomorphism of the even dimensional nilmanifolds constructed in Theorem \ref{fpf-maps-on-nilmanifolds} can be isotoped to a fixed point free diffeomorphism but the isotopy does not respect the Hamiltonian structure so that the resulting fixed point free diffeomorphism is {\it not} a Hamiltonian symplectomorphism.
\end{remark}

\begin{remark}
It is well-known (see e.g. \cite{charlap} or \cite{vasq}) that there is a one-to-one correspondence between compact flat manifolds and finitely generated Bieberbach groups, that is, compact Riemannian manifolds with zero sectional curvatures are precisely the finite dimensional aspherical manifolds with finitely generated torsion-free virtually abelian fundamental groups. It is natural to ask whether one can obtain a result similar to Theorem \ref{fpf-maps-on-nilmanifolds} for flat manifolds. One of the main differences between nilmanifolds and infra-nilmanifolds (or flat manifolds) is that $N(f)=0$ iff $R(f)=\infty$ for selfmaps on nilmanifolds whereas there are maps on flat manifolds where $R(f)=\infty$ but $N(f)>0$ (see Example 5.3 of \cite{daci-peter2}). Even for the Klein bottle $K$, one can find a homeomorphism with $N(f)>0$ and $R(f)=\infty$. Thus, the groups $\pi_1(K) \times \mathbb Z^n$ do not provide examples of flat manifolds on which all homeomorphisms have zero Nielsen numbers.
\end{remark}

A natural extension of our results for nilmanifolds is that of an infra-nilmanifold whose fundamental group is a finitely generated torsion-free virtually nilpotent group or equivalently an almost Bieberbach group. From the point of view of Nielsen fixed point theory, D. Anosov in \cite{anosov} (or \cite{ed-suf}) already pointed out that the equality $N(f)=|L(f)|$ does not hold for selfmaps of infra-nilmanifolds. Furthermore, Kwasik and Lee \cite{kwasik-lee} constructed infra-nilmanifolds and affine Anosov diffeomorphisms $f$ for which $|L(f)|\ne N(f)$ in every even dimension $n\ge 4$. On the other hand, for a selfmap $f$ of compact solvmanifold, $R(f)<\infty$ implies that $N(f)=R(f)$ \cite{daci-peter2}. Thus, one would like to study the $R_\infty$ property for such groups, although we cannot always obtain results similar to those in the previous two sections.

\section{$\mathcal C$-nilpotent groups}

Next we turn to another generalization of nilpotent groups. Nilpotent spaces have proven to be useful in homotopy theory and in algebraic topology in general. The concept of a Serre class of abelian groups has been generalized. One such generalization is that of a $\mathcal C$-nilpotent class, first introduced in \cite{Daci83}. A family $\mathcal C$ of groups is a {\it class of groups} if given a short exact sequence of groups
$$
1\to A \to B \to C \to 1
$$
then $A,C\in {\mathcal C}$ if and only if $B\in {\mathcal C}$. Given a class ${\mathcal C}$, a group $G$ is said to be ${\mathcal C}$-nilpotent if $\Gamma_n(G)\in {\mathcal C}$ for some positive integer $n$ where $\Gamma_i(G)$ is the $i$-th commutator in the lower central series of $G$. More generally, a group $\pi$ acts on a group $G$ $\mathcal C$-nilpotently if $\Gamma_n^{\pi}(G)\in {\mathcal C}$ for some positive integer $n$ where $\Gamma_n^{\pi}(G)$ is the smallest $\pi$-subgroup that contains $[G,\Gamma_{n-1}^{\pi}(G)]$ and the set $\{(\alpha\cdot g)g^{-1}|\alpha \in \pi, g\in \Gamma_{n-1}^{\pi}(G)\}$. Thus, a space $X$ is said to be $\mathcal C$-nilpotent if $\pi_1(X)$ is $\mathcal C$-nilpotent and the action of $\pi_1(X)$ on $\pi_k(X)$ is $\mathcal C$-nilpotent. From now on, we let ${\mathcal C}$ be the class of finite groups.

Let $G$ be a finitely generated group. For any $\varphi \in {\rm Aut}(G)$ and any positive integer $n$, we have a commutative diagram
\begin{equation}
\begin{CD}
    0 @>>> \Gamma_n(G) @>>> G           @>>>      G/\Gamma_n(G)        @>>> 1   \\
    @.     @V{\varphi'}VV      @V{\varphi}VV          @V{\bar \varphi}VV          @.  \\
    0 @>>> \Gamma_n(G) @>>> G           @>>>      G/\Gamma_n(G)        @>>> 1
\end{CD}
\end{equation}

It follows easily from Lemma \ref{R-facts} that if $G/\Gamma_n(G)$ has the $R_{\infty}$ property then so does $G$. Since $N_n=G/\Gamma_n(G)$ is nilpotent, the torsion elements form a subgroup $T_n$ and $T_n$ is a characteristic subgroup of $G/\Gamma_n(G)$. It follows that if $N_n/T_n$ has the $R_{\infty}$ property then $N_n$ also has the $R_{\infty}$ property.

Thus, we have the following useful result that can be used to construct more examples of groups with the $R_{\infty}$ property.

\begin{proposition}\label{C-nil}
Let $G$ be a finitely generated group and $T_k$ be the torsion subgroup of the nilpotent group $N_k=G/\Gamma_k(G)$. If for some positive integer $n$ the group $N_n/T_n$ has the $R_{\infty}$ property then $G$ also has the $R_{\infty}$ property.
\end{proposition}

We now combine our results from section 4 to enlarge the family of groups with the $R_{\infty}$ property and also the family of compact manifolds which have fundamental groups with the $R_{\infty}$ property.

\begin{example}\label{Poincare-nil}
Let $\Sigma^3$ be the Poincar\'e $3$-sphere with $\pi_1(\Sigma^3)\cong Icos$, the binary icosahedral group of order $120$. Suppose $M$ is the $4$-dimensional $K(\pi,1)$ nilmanifold where $\pi$ is the torsion-free nilpotent group of nilpotency class $3$ as in Example \ref{ex3}. Take $X=\Sigma^3\times M$. Now, $G=\pi_1(X)$ is $\mathcal C$-nilpotent, $\Gamma_4(G)=Icos$, and $G/\Gamma_4(G)\cong \pi$. Since $\pi$ has the $R_{\infty}$ property (and the torsion subgroup of $\pi$ is trivial), it follows from Proposition \ref{C-nil} that $G$ has the $R_{\infty}$ property. More generally, one can obtain similar examples by replacing $Icos$ with any finite perfect group.
\end{example}

We can show more for homeomorphisms of the space $X$ in the last example although it is not clear if $X$ is a Jiang-type space. In fact, we obtain the following result, similar to Theorem \ref{fpf-maps-on-nilmanifolds}.

\begin{theorem}\label{fpf-maps-on-C-nilpotent-spaces}
For any $n\ge 7$, there exists a compact $\mathcal C$-nilpotent manifold $M$ with $\dim M=n$ such that every homeomorphism $f:M\to M$ is isotopic to a fixed point free homeomorphism.
\end{theorem}
\begin{proof} Let $n\ge 7$ and $M=\Sigma^3 \times N^{n-3}$ where $N^k$ is the compact nilmanifold with fundamental group of Hirsch length $k$ as constructed in Example \ref{ex3}.
We will show that for any  $f:M \to M$, $f$ is homotopic to a fiber-preserving map. First, choose basepoints $x_0\in \Sigma^3$, $y_0\in N=N^{n-3}$ and write $f(x,y)=(g(x,y),h(x,y))$. Denote by $\bar f:N\to N$ the restriction of $h$ on $\{x_0\}\times N$, that is, $\bar f(y)=h(x_0,y)$. Next, we will show that $\bar f\circ p_2$ and $h$ are homotopic where $p_2:\Sigma^3\times N\to N$ is the projection on the second factor. Write $A=\Sigma^3\times N\times \{0,1\} \cup (\Sigma^3\vee M)\times [0,1]$. The pair $(M,A)$ is a relative $CW$ complex. Since $\pi_1(\Sigma^3)\cong Icos$ is finite, $N$ is aspherical and $\pi=\pi_1(N)$ is torsion-free, the maps $h$ and $\bar f\circ p_2$ restricted to $\Sigma^3$ are null-homotopic and they are homotopic on $\{x_0\}\times N$. It follows that $h$ and $\bar f\circ p_2$ coincide, up to homotopy, on the subspace $\Sigma^3 \vee N \subset \Sigma^3 \times N$. Let $F$ be the homotopy from $h|_{\Sigma^3\vee N}$ to $(\bar f\circ p_2)|_{\Sigma^3\vee N}$ and $\widehat F:A\to N$ where $\widehat F|_{M\times \{0\}}=h, \widehat F|_{M\times \{1\}}=\bar f\circ p_2$ and $\widehat F|_{A\times [0,1]}=F$. Since the 2-skeleton of $M=\Sigma^3\times N$ already lies inside $\Sigma^3\vee N\subset A$, there is a sequence of obstructions $c^i(\widehat F)\in H^i(M,A;\pi_{i-1}(N))$ to extending $\widehat F$ to $M$ where $i>2$. These obstructions are trivial because $\pi_k(N)=0$ for $k>1$ since $N$ is aspherical. Now, up to homotopy, the map $h$ is of the form $\bar f\circ p_2$ thus we have the following commutative diagram
\begin{equation}\label{fiber}
\begin{CD}
    M    @>{f}>>  M  \\
    @V{p_2}VV  @VV{p_2}V   \\
    N    @>{\bar f}>>  N
\end{CD}
\end{equation}
Now, suppose $f$ is a homeomorphism and \eqref{fiber} holds.
If $f$ does not have any fixed points there is nothing to prove. Suppose $f$ does have a fixed point then \eqref{fiber} yields a commutative diagram of groups similar to the diagram as in Lemma \ref{R-facts}. Thus, $\bar f$ induces an automorphism $\bar \varphi$ which has $R(\bar \varphi)=\infty$ since $\pi$ has the $R_{\infty}$ property. It follows that $N(\bar f)=0$. Since the fibration $p_2$ is trivial, it follows that $N(f)=N(f')\cdot N(\bar f)=0$ where $f'$ is the restriction of $f$ on the fiber. Since $n\ge 7$, Kelly's Isotopy Wecken Theorem \cite{kelly} implies that $f$ is isotopic to a fixed point free homeomorphism.
\end{proof}

\begin{remark}
Although the space $M$ in Theorem \ref{fpf-maps-on-C-nilpotent-spaces} is a ${\mathcal C}$-nilpotent space, its fundamental group has a center of infinite index so we cannot conclude from \cite{daci-peter3} that $M$ is of Jiang-type. We should point out that one can construct similar manifolds $M$ with the desired property by replacing $\Sigma^3$ with any compact manifold $L$ with finite $\pi_1(L)$ since the same argument also shows that every map $f:L\times N \to L\times N$ is fiber-preserving up to homotopy. Note that if $\dim L\le 2$, $\pi_1(L)$ is necessarily abelian so that $\pi_1(M)$ would be nilpotent possibly with torsion.
\end{remark}

\begin{remark}
It is shown in \cite{daci-peter3} that if $H$ is a finite subgroup of a compact connected Lie group $G$, then the coset space $X=G/H$ is $\mathcal C$-nilpotent and its fundamental group has a center of finite index and hence is of Jiang-type. However, these spaces do not provide examples like Example \ref{Poincare-nil} because $[\pi_1(X),\pi_1(X)]$ is also finite. It follows that $\pi_1(X)/[\pi_1(X),\pi_1(X)]=\pi_1(X)^{Ab}$ is finitely generated abelian and therefore does not have the $R_{\infty}$ property and hence Proposition \ref{C-nil} is not applicable.
\end{remark}


\begin{thebibliography}{99}

\bibitem{anosov} D. Anosov, The Nielsen number of maps of
nilmanifolds,
{\em Russian Math. Surveys} {\bf 40} (1985), 149--150.

\bibitem{bryant-papistas} R.M. Bryant and A. Papistas, Automorphism groups of nilpotent groups, {\em
Bull. London Math. Soc.} {\bf 21} (1989), 459--462.

\bibitem{charlap} L. Charlap, ``Bieberbach Groups and Flat Manifolds," Springer,
New York, 1986.

\bibitem{dyer} J. Dyer, A nilpotent Lie algebra with nilpotent automorphism group, {\em
Bull. Amer. Math. Soc.} {\bf 76} (1970), 52--60.

\bibitem{ed-suf} E. Fadell and S. Husseini, On a theorem of Anosov on Nielsen numbers for Nilmanifolds,{Nonliner functional analizes and its applicatios} (Maratea, 1985) NATO, Advanced Science Institutes Series C: Mathematical and  Physical  Science 173(Kluwer, Dordrecht, 1986), 47-53.

\bibitem{fel:1} A. L. Fel'shtyn, The Reidemeister number of any automorphism
of a Gromov hyperbolic group is infinite,
{\em Zapiski Nauchnych Seminarov POMI} {\bf 279} (2001), 229--241.

\bibitem{fel-daci} A. L. Fel'shtyn and D. Gon\c calves, Reidemeister numbers of any automorphism of Baumslag-Solitar groups is infinite, in: Geometry and Dynamics of Groups and Spaces, Progress in Mathematics, v.265 (2008), 286--306.

\bibitem{fgw} A. L. Fel'shtyn, D. Gon\c calves and P. Wong, Twisted conjugacy
classes for polyfree groups, preprint 2007.

\bibitem{fel-hill} A. L. Fel'shtyn and R. Hill, The Reidemeister zeta function
with applications to Nielsen theory and a connection with Reidemeister
torsion,  {\em K-theory}
{\bf 8} no.4 (1994), 367--393.

\bibitem{felshtyn-hill-wong} A. L. Fel'shtyn, R. Hill and P. Wong, Reidemeister numbers
of equivariant maps, {\em Top. Appl.} {\bf 67} (1995), 119--131.

\bibitem{felshtyn-leonov-troitsky} A. L. Fel'shtyn, Y. Leonov and E. Troitsky, Reidemeister numbers of saturated weakly branch groups, arXiv:math/0606725, preprint 2006.

\bibitem{felshtyn-troitsky} A. L. Fel'shtyn and E. Troitsky, Twisted Burnside-Frobenius theory for discrete groups, {\em J. reine angew. Math.}, to appear.

\bibitem{go:nil2} D. L. Gon\c calves, The coincidence Reidemeister classes of
maps on nilmanifolds, {\em Topol. Methods Nonlinear Anal.} {\bf 12}
no.2 (1998), 375--386.

\bibitem{Daci83} D. Gon\c calves, Generalized classes of groups, $C$-nilpotent spaces and "the Hurewicz theorem", {\em Math. Scand.}  {\bf 53}  (1983), 39--61.

\bibitem{daci-peter} D. Gon\c calves and P. Wong, Twisted conjugacy classes in exponential growth groups, {\em Bull. London Math. Soc.} {\bf 35} (2003), 261--268.

\bibitem{daci-peter2} D. Gon\c calves and P. Wong, Homogeneous spaces in coincidence theory. II,  {\em Forum Math.}  {\bf 17} (2005), 297--313.

\bibitem{daci-peter3} D. Gon\c calves and P. Wong, Homogeneous spaces in coincidence
theory, in: P. Schweitzer, ed., Conference proceedings of the Tenth
Brazilian Topology Meeting (S\~ao Carlos 1996), {\em Mat. Contemp.} {\bf 13}
(1997), 143--158.

\bibitem{daci-peter4} D. Gon\c calves and P. Wong, Twisted conjugacy classes in wreath products,  {\em Internat. J. Alg. Comput.}  {\bf 16}  (2006), 875--886.

\bibitem{jww} B. Jiang, S. Wang, and Y. Wu,
Homeomorphisms of $3$-manifolds and the realization of Nielsen number,
{\em Comm. Anal. Geom.} {\bf 9} (2001), 825--877.

\bibitem{jiang-guo} B. Jiang and J. Guo, Fixed points of surface diffeomorphisms, {\em Pacific J. Math.} {\bf 160} (1993), 67--89.

\bibitem{kelly} M. Kelly, The Nielsen number as an isotopy invariant,
{\em Topology Appl.} {\bf 62} (1995), no. 2, 127--143.


\bibitem{kwasik-lee} S. Kwasik and K. B. Lee, The Nielsen numbers of homotopically periodic maps of  infranilmanifolds, {\em J. London Math. Soc.} {\bf 38}  (1988), 544--554.

\bibitem{levitt} G. Levitt, On the automorphism group of generalized Baumslag-Solitar groups, {\em Geometry and Topology} {\bf 11} (2007), 473--515.

\bibitem{levitt-lustig} G. Levitt and M. Lustig, Most automorphisms of a
hyperbolic group have very simple dynamics, {\em Ann. Scient. \'Ec. Norm. Sup.}
{\bf 33} (2000), 507--517.

\bibitem{kms} W. Magnus, A. Karrass and D. Solitar, Combinatorial Group
Theory. New York: Interscience Publishers, (1966).

\bibitem{MO} C. McCord and J. Oprea, Rational Lusternik-Schnirelmann category and the Arnold conjecture for nilmanifolds, {\em Topology} {\bf 32} (1993), 701--717.

\bibitem{R} Y. Rudyak, On analytical applications of stable homotopy (the Arnold conjecture, critical points) {\em Math. Zeitschrift} {\bf 230} (1999), 659--672.

\bibitem{TW1}
J. Taback and P. Wong, Twisted conjugacy and quasi-isometry invariance for generalized solvable Baumslag-Solitar groups, {\em Journal London Mathematical Society (2)} {\bf 75} (2007), 705-717.

\bibitem{TW2}
J. Taback and P. Wong, A note on twisted conjugacy and generalized Baumslag-Solitar groups, preprint, 2006.

\bibitem{vasq} A. T. Vasquez, Flat Riemannian manifolds, {\em J. Diff. Geom.} {\bf 524} (1970), 367--382.

\bibitem{wong} P. Wong, Reidemeister number, Hirsch rank,
coincidences on polycyclic groups and solvmanifolds,
{\em J. reine angew. Math.} {\bf 524} (2000), 185--204.


\end{thebibliography}
\end{document}